\UseRawInputEncoding
\documentclass{article}
\usepackage{amsfonts}
\usepackage{amscd}
\usepackage{latexsym,bm,amssymb}
\usepackage{amsmath}
\usepackage{amsthm}
\usepackage{indentfirst}
\usepackage[sort&compress,numbers]{natbib}
\usepackage{color}
\usepackage{multirow}
\usepackage{graphicx}
\usepackage[noend]{algpseudocode}
\usepackage{algorithmicx,algorithm}
\usepackage{extarrows}
\usepackage[perpage,symbol*]{footmisc}
\usepackage{footmisc}
\usepackage{subfigure}
\newtheorem{definition}{Definition}[section]
\newtheorem{theorem}{Theorem}[section]
\newtheorem{lemma}{Lemma}[section]
\newtheorem{remark}{Remark}[section]

\begin{document}
\newcommand{\re}[1]{\begin{rema} {\rm{#1}} \end{rema}}
\newcommand{\de}[1]{\begin{Def} {\rm{#1}} \end{Def}}
\newcommand{\ex}[1]{\begin{exm} {\rm{#1}} \end{exm}}
\newcommand{\norm}[1]{\parallel #1\parallel}
\newcommand{\seq}[1]{\langle #1\rangle}
\newcommand{\sty}{\displaystyle}
\newcommand{\ra}{\rightarrow}
\newcommand{\Ra}{\Rightarrow}
\newcommand{\n}{ {\mathcal N} }
\def\QED{\hfill $\Box$\smallskip}

\title{{\Large {\bf An inexact primal-dual method with correction step for a saddle point problem in image debluring}}
\thanks{This work is partially supported by the National Natural Science Foundation of
China (Nos. 11771350 and 12101100), Basic and Advanced Research Project of CQ CSTC (No.  cstc2020jcyj-msxmX0738) }
\author{Changjie Fang,\thanks{Corresponding author, E-mail address: fangcj@cqupt.edu.cn} \, Liliang Hu,\,Shenglan Chen\\
\small\it Key Lab of Intelligent Analysis and Decision on Complex Systems, \\\small\it Chongqing University
of Posts and Telecommunications, Chongqing 400065, China\\\small\it  School of Science, Chongqing University
of Posts and Telecommunications,\\
\small\it Chongqing
400065,
 China}}
\date{}
\maketitle
\vspace*{-9mm}
\begin{center}
\begin{minipage}{5.5in}
{\bf Abstract.}\quad  In this paper, we present an inexact primal-dual method with correction step for a saddle point problem by introducing the notations of inexact extended proximal  operators with symmetric positive definite matrix $D$. Relaxing  requirement on  primal-dual step sizes, we prove the convergence of the proposed method. We also establish the $O(1/N)$ convergence rate of our method in the ergodic sense. Moreover, we apply our method to solve TV-L$_1$ image deblurring problems. Numerical simulation results illustrate the efficiency of our method.
\\ \ \\
{\bf Key Words.}\quad primal-dual method, inexact extended proximal operators, convergence rate, prediction and correction, image deblurring
\\ \ \\
\end{minipage}
\end{center}
\section{Introduction}
Let  $X$ and $Y$ be two finite-dimensional real vector spaces equipped with  an inner product $\langle\cdot,\cdot\rangle$ and a norm $\|\cdot\|=\sqrt{\langle\cdot,\cdot\rangle}$.
In this paper,we consider the following saddle point problem:
\begin{align}\label{1.1}
\min_{x \in X} \max_{y \in Y} L(x,y)=f(x)+\langle Ax,y\rangle-g(y)
\end{align}
where $A$ is a bounded linear mapping denoting some imaging operator, $f: X\rightarrow(-\infty,+\infty]$  and  $g: Y\rightarrow(-\infty,+\infty]$ are proper lower semicontinuous (l.s.c) convex functions.

Recall that $(x^*,y^*)$ is called the saddle point of \eqref{1.1}, if
\begin{align}\label{1.2}
 L(x^*,y)\leq L(x^*,y^*)\leq L(x,y^*),\forall x\in X,y\in Y
\end{align}
Now  we consider the primal problem
 \begin{align}\label{1.2a}
 \min_{x\in X}f(x)+h(Ax)
\end{align}
together with its dual problem
\begin{align}\label{1.2b}
 \max_{y\in Y}-f^*(-A^*y)-h^*(y)
\end{align}
where $h^*$ denotes the Legendre-Fenchel conjugate of a convex l.s.c. function $h$, $A^*$  denotes the adjoint of the bounded linear operator $A$.

If a primal-dual solution pair $(x^*,y^*)$ of \eqref{1.2a} and \eqref{1.2b} exists, i.e.,
\begin{align*}
 0\in\partial f(x^*)+A^*y^*,\,\, 0\in\partial h(Ax^*)-y^*,
\end{align*}
then the problem \eqref{1.2a} is equivalent to the following saddle-point formulation:
\begin{align}\label{1.2c}
  \min_{x\in X}\max_{y\in Y}f(x)+\langle Ax,y\rangle-h^*(y).
\end{align}
Hence, Problem \eqref{1.2c} is a special case of Problem \eqref{1.1}.

It is well known that many application problems  can be
formulated as the saddle point problem \eqref{1.1} such as image restoration, magnetic resonance imaging and computer vision; see, for example, \cite{PCBC2009,ROF1992,Val2014,ZC2008}.

Two of the most popular approaches are first-order primal-dual methods \cite{CP2011,EZC2010}, in particular the Primal-Dual Hybrid Gradient (PDHG) method \cite{ZC2008}, and Alternating Direction Method of Multipliers (ADMM) method \cite{BPCPE2010,HTY2012}. For PDHG method, both a primal and a dual variable are updated in each iteration and thus some difficulties  that arise when working only on the primal or dual variable can be avoided. In  ADMM method, separating the minimization over the two primal variables into two steps is precisely what allows for decomposition when $f$ or $g$, or both, are separable. In \cite{HYY2014}, it was showed that  PDHG is not necessarily convergent even when the step sizes are fixed as tiny constants. In \cite{EZC2010}, PDHG was interpreted as projected averaged gradient method, and its convergence was studied
 by imposing additional restrictions ensuring that the step sizes $\lambda$ and $\tau$ are small. In \cite{CP2011}, a primal-dual method with inertial step $\theta$($\theta\in [0,1]$) was proposed (denoted by PDI) with convergence rate $O(1/N)$ in terms of primal-dual gap, and for $\theta=1$, the convergence of PDI was proved with the
requirement on step sizes $\tau\lambda<1/\|A^TA\|$. In \cite{HY2012}, some prediction-correction contraction methods  were presented in which the convergence was guaranteed with relaxed step sizes satisfying $\lambda\tau<4/(1+\theta)^2\|A^TA\|^2$ for $\theta\in(-1,1]$. Further, a primal-dual method (named PDL) in the prediction-correction fashion was proposed in \cite{Xie2017} and the pairwise primal-dual stepsizes $\lambda_i$ and $\tau_i$ were relaxed to $\lambda_i\tau_i<1(i=1,2)$.

As first-order methods, however, they are sensitive to problem conditions, and hence might be performed up to a certain precision, for example, due to the application of a proximal operator lacking a closed-form solution. This problem may arise from examples studied in, for example, \cite{CCS2010,EB2016,FP2011,MGC2011}. An absolute error criterion was
adopted in \cite{EB1992}, where the subproblem errors are controlled by a summable  sequence of error tolerances. To simplify the choice of the sequences, a relative error criterion was  introduced in \cite{NWY2011}, where the corresponding parameters are required to be square summable.  \cite{RC2020}  introduced four different types of inexact proxima, where all the controlled errors were required to be summable. In \cite{LXY2020}, the inexact  preconditioned PDHG method was studied by the selection of appropriate preconditioners and the introduction of bounded relative error of the subproblem, where  convergence was established in case the error was neither summable nor square summable.

Motivated by the research works \cite{RC2020,Xie2017}, in this paper, we introduce three different types of  inexact extended proxima closely related to the  extended  proximal operator with the matrix $D$(\cite{CR2013}). Applying these notions, we propose an inexact primal-dual method with correction step for solving the saddle point problem. Under some mild conditions, the convergence of the proposed method is proved, in which we relax  requirement on pairwise primal-dual stepsize , for example, compared with that in \cite{Xie2017}. In \cite{Xie2017}, primal-dual stepsizes $\lambda_i$ and $\tau_i(i=1,2)$  are required to satisfy $\lambda_i\tau_i<1$. In our method, the sequences $\{\lambda_k\}$ and $\{\tau_k\}$ are  nondecreasing and bounded satisfying $R-\tau_k\lambda_kS^{-1}\succ 0 $, where $R$ and $S$ are  symmetric positive definite matrices. We also establish the $O(1/N)$ convergence rate in the  ergodic sense. At the same time, we  establish the convergence rates in case  error tolerances $\{\delta_k\}$ and $\{\varepsilon_k\}$ are required to  decrease like $O(1/k^{2\alpha+1})$ for some $\alpha>0$; see Theorem \ref{th3.3}. In the numerical experiments part, we investigate the applications of our method in TV-L$_1$ image deblurring. Firstly, we show that the type-2 approximation of the extended proximal point can be computed by approximately minimizing  duality gap; see \eqref{4.5a} in Section \ref{sec:ex}. Further, the duality gap is used as the  stopping criterion of inner loop, i.e., the second subproblem; see \eqref{4.6}. In addition, we discuss the sensitivity of parameters in Algorithm 1. Finally, we show the efficiency of our method in image deblurring compared with some existing  methods, for example, \cite{CP2011,RC2020,Xie2017}.

The rest of this paper is organized as follows. In Section \ref{sec:pr}, we introduce the concepts of  inexact extended proximal  operators and present some auxiliary  lemmas. In Section \ref{sec:alg},  we describe our method and prove the convergence of our method. At the same time, we also analyze the convergence rate. Numerical experiment results are reported in Section \ref{sec:ex}. Some conclusions are presented in Section \ref{sec:con}.
\section{Preliminaries}\label{sec:pr}
\setcounter{equation}{0}
In this section, we shall introduce some definitions. Suppose that $h$ be a convex function in $\mathbb{R}^n$, $D \in \mathbb{R}^{n\times n}$ a  symmetric positive definite matrix and $\tau>0$. For any $D\succ0$ and  given $y \in \mathbb{R}^n$ , denote
\begin{align}\label{2.1}
G_y(x):=h(x)+\frac{1}{2\tau}\|x-y\|_D^2,\,\forall\,\,x\in \mathbb{R}^n,
\end{align}
and define the extended proximal operator of $h$ as
\begin{align}\label{2.2}
\textrm{Prox}_{\tau h}^D(y):=\overline{z}=\arg\min_{x \in X}G_y(x),
\end{align}
where $\|x\|_D^2= \langle x,Dx \rangle$ and $D^{-1}$ denotes the inverse of $D$. Because $D$ is symmetric positive definite, $\textrm{Prox}_{\tau h}^D(y)$ is unique(see Lemma \ref{le2.4}).
\begin{definition}\label{d2.1} Let $\varepsilon\geq0$. $z\in X$ is said to be a type-0 approximation of the extended proximal point $\emph{Prox}_{\tau h}^D(y)$ with precision $\varepsilon$ if
\begin{align*}
z\approx_0^{\varepsilon} \emph{Prox}_{\tau h}^D(y)\overset{def}\Longleftrightarrow\|z-\bar{z}\|_D\leq \sqrt{2\tau\varepsilon}
\end{align*}
\end{definition}
Next we recall the definition of $\varepsilon$-subdifferential of $h$ at $z$, denoted by $\partial_{\varepsilon}h(z)$:
\begin{align*}
\partial_{\varepsilon}h(z)=\{p\in X|h(x)\geq h(z)+\langle p,x-z\rangle-\varepsilon ,\forall x \in X\}.
\end{align*}
In the following, we give the definition of $\varepsilon$-subdifferential of $G_y$ at $z$, denoted by  $\partial_{\varepsilon}G_y(z)$:
\begin{align*}
\partial_{\varepsilon}G_y(z):=\{p\in X|G_y(x)\geq G_y(z)+\langle p,x-z\rangle-\varepsilon ,\forall x \in X\}.
\end{align*}
\begin{definition}\label{d2.2}Let $\varepsilon\geq0$. $z\in X$ is said to be a type-1 approximation of the extended proximal point $\emph{Prox}_{\tau h}^D(y)$ with precision $\varepsilon$ if
\begin{align*}
z\approx_1^{\varepsilon} \emph{Prox}_{\tau h}^D(y)\overset{def}\Longleftrightarrow0 \in \partial_{\varepsilon}G_y(z)
\end{align*}
\end{definition}
\begin{definition}\label{d2.3} Let $\varepsilon\geq0$. $z\in X$ is said to be a type-2 approximation of the extended proximal point $\emph{Prox}_{\tau h}^D(y)$ with precision $\varepsilon$ if
\begin{align*}
z\approx_2^{\varepsilon} \emph{Prox}_{\tau h}^D(y)\overset{def}\Longleftrightarrow\frac{1}{\tau}D(y-z)\in \partial_{\varepsilon}h(z)
\end{align*}
\end{definition}
\begin{remark}\label{re2.1}
If $D=I$, where $I$ is the identity matrix, then  the inexact extended proximal operators in Definitions \ref{d2.1}-\ref{d2.3} will reduce into the inexact  proxima, for example, introduced in \cite{RC2020}, respectively. Thus, Definitions \ref{d2.1}-\ref{d2.3} are the generalization of the correcponding definitions in \cite{RC2020}.
\end{remark}
According to the above definitions, we have the following lemmas.
\begin{lemma}\label{le2.1}Suppose
$z\approx_1^{\varepsilon}\mathop{\arg\min}_{x \in X}\{h(x)+\frac{1}{2\tau}\|x-y\|_D^2\}$, then $z\in$ dom $h$ and $z\approx_0^{\varepsilon} \arg\min_X\{h(x)+\frac{1}{2\tau}\|x-y\|_D^2\}$.
\end{lemma}
\begin{proof}According to Definition \ref{d2.2} and the definition of  $\partial_{\varepsilon}G_y$, we have
\begin{align}\label{g01}
G_y(x)\geq G_y(z)-\varepsilon,\quad \forall x \in X.
\end{align}
Setting $x=\overline{z}$ in \eqref{g01} and using \eqref{2.2}, from \eqref{g01} we have
\begin{align} \label{g02}
\|\overline{z}-z\|^2_D \leq 2\tau\varepsilon+2\tau[h(\overline{z})-h(z)+\langle \frac{1}{\tau}D(\overline{z}-y),\overline{z}-z\rangle]，
\end{align}
which implies that $z \in dom f$. According to the optimality condition of \eqref{2.2}, we have
\begin{align}\label{h01}
h(x) \geq h(\overline{z})+\langle \frac{1}{\tau}D(y-\overline{z}),x-\overline{z}\rangle,\quad \forall x \in X.
\end{align}
Setting $x=z$ in \eqref{h01} and substituting the resulting inequality into \eqref{g02}, we get
\begin{align*}
||\overline{z}-z||_D \leq \sqrt{2\tau\varepsilon}.
\end{align*}
In view of  Definite \ref{d2.1}, we obtain the conclusion.
\end{proof}
\begin{lemma}\label{le2.2}Suppose that$z\approx_1^{\varepsilon} \arg\min_{x \in X}\{h(x)+\frac{1}{2\tau}\|x-y\|_D^2\}$, then there exists $r\in X$ with$\|r\|_D\leq \sqrt{2\tau\varepsilon}$ such that
\begin{align*}
\frac{1}{\tau}D(y-z-r)\in \partial_{\varepsilon}h(z).
\end{align*}
\end{lemma}
\begin{proof} According to Definition \ref{d2.2}, from \eqref{2.1} we have
\begin{align*}
h(x) \geq h(z)+\langle x-z,\frac{1}{\tau}D(y-z-\frac{x-z}{2})\rangle-\varepsilon,\quad \forall x \in X.
\end{align*}
Set $r=\frac{x-z}{2}$. By the definition of the $\varepsilon$-subdifferential of $h$, the conclusion holds.
\end{proof}
\begin{lemma}\label{le2.3}Suppose $z\approx_2^{\varepsilon} \arg\min_{x\in X}\{h(x)+\frac{1}{2\tau}\|x-y\|_D^2\}$, then $z\approx_1^{\varepsilon}\arg\min_{x \in X}\{h(x)+\frac{1}{2\tau}\|x-y\|_D^2\}$
\end{lemma}
\begin{proof} According to Definition \ref{d2.3}, from \eqref{2.1} we have
\begin{align*}
G_y(x)&=h(x)+\frac{1}{2\tau}\|x-y\|_D^2\\
         &\geq h(z)+\langle\frac{1}{\tau}D(y-z),x-z\rangle-\varepsilon+\frac{1}{2\tau}\|x-y\|_D^2\\
         &=h(z)+\frac{1}{2\tau}\|z-y\|_D^2+\frac{1}{2\tau}\|x-z\|_D^2-\varepsilon \\
         &=G_y(z)-\varepsilon+\frac{1}{2\tau}\|x-z\|_D^2\\
         &\geq G_y(z)-\varepsilon,
\end{align*}
where the first inequality follows from the definition of $\partial_{\varepsilon}h(z)$ and  the second equality follows from following identity equality
\begin{align}\label{equ}
\|a-b\|_D^2=\|a-c\|_D^2+\|b-c\|_D^2-2\langle D(a-c),b-c\rangle.
\end{align}
Hence, $0\in\partial_{\varepsilon}G_y(z)$. By Definition \ref{d2.2}, the conclusion holds.
\end{proof}
The following lemma illustrates that the extended proximal operator \eqref{2.2} is well-defined.
\begin{lemma}\label{le2.4}
\begin{align*}
\|\emph{Prox}_{\tau h}^D(y_1)-\emph{Prox}_{\tau h}^D(y_2)\|_D\leq \|y_1-y_2\|_D.
\end{align*}
\end{lemma}
\begin{proof}Let $z_1=\textrm{Prox}_{\tau h}^D(y_1)$ and $z_2=\textrm{Prox}_{\tau h}^D(y_2)$. Since $z_1=\textrm{Prox}_{\tau h}^D(y_1)$, the optimality condition implies that $\frac{1}{\tau}D(y_1-z_1) \in \partial h(z_1)$, and hence
\begin{align}\label{p01}
h(x)\geq h(z_1)+\langle \frac{1}{\tau}D(y_1-z_1),x-z_1\rangle,\quad \forall x \in X.
\end{align}
Set $x=z_2$ in \eqref{p01} and get
\begin{align}\label{2.4}
h(z_2)\geq h(z_1)+\langle \frac{1}{\tau}D(y_1-z_1),z_2-z_1\rangle.
\end{align}
Similarly, we have
\begin{align}\label{2.5}
h(z_1)\geq h(z_2)+\langle \frac{1}{\tau}D(y_2-z_2),z_1-z_2\rangle.
\end{align}
Adding \eqref{2.4} and \eqref{2.5}, and by a simple manipulation, we obtain
\begin{align*}
\|z_1-z_2\|_D^2 \leq \langle D(y_2-y_1),z_2-z_1\rangle \leq\|y_2-y_1\|_D\|z_2-z_1\|_D,
\end{align*}
i.e.,
\begin{align*}
\|\textrm{Prox}_{\tau h}^D(y_1)-\textrm{Prox}_{\tau h}^D(y_2)\|_D \leq\|y_1-y_2\|_D.
\end{align*}
\end{proof}
The following lemma is crucial in proving the convergence of Algorithm 1.
\begin{lemma}\label{le2.5}Suppose that $g:X\mapsto\bar{R}$ is a convex function. For given  $z_0,u,v \in X$ and $\tau,\varepsilon>0$. Let
\begin{align}
z_1&\approx_2^{\varepsilon} \arg\min_{z\in X}\{g(z)+\frac{1}{2\tau}\|z-(z_0-D^{-1}u)\|_D^2\}, \label{2.6}\\
z_2&\approx_1^{\varepsilon} \arg\min_{z\in X}\{g(z)+\frac{1}{2\tau}\|z-(z_0-D^{-1}v)\|_D^2\} , \label{2.7}
\end{align}
then
\begin{itemize}
\item[(i)]
\begin{align} \label{2.9}
0\leq \|z_1-z_2\|_D\leq \frac{1}{2}( \sqrt{2\tau\varepsilon}+\|u-v\|_{D^{-1}}+\sqrt{\|u-v\|_{D^{-1}}^2+10\tau\varepsilon+2\sqrt{2\tau\varepsilon}\|u-v\|_{D^{-1}}}).
\end{align}
\item[(ii)] $\forall z\in X$,
\begin{align}
g(z_1)-g(z)+&\langle z_1-z,\frac{1}{\tau}v \rangle\leq \frac{1}{2\tau}(\|z-z_0\|_D^2-\|z-z_2\|_D^2) \notag \\
&+\frac{1}{2\tau}(\|D^{-1}(u-v)\|_D^2-\|z_0-z_1\|_D^2)+2\varepsilon+\sqrt{\frac{2\varepsilon}{\tau}}\|z-z_2\|_D.   \label{2.8}
\end{align}
\end{itemize}
\end{lemma}

\begin{proof}(i) In view of Definition \ref{d2.3} and \eqref{2.6}, we have
\begin{align}
\frac{1}{\tau}D((z_0-D^{-1}u)-z_1) \in \partial_{\varepsilon}g(z_1). \label{2.10}
\end{align}
By Lemma \ref{le2.2}, there exists $r\in X$ with$||r||_D\leq \sqrt{2\tau\varepsilon}$ such that
\begin{align}
\frac{1}{\tau}D((z_0-D^{-1}v)-z_2-r)\in \partial_{\varepsilon}g(z_2). \label{2.11}
\end{align}
 From the definition of $\varepsilon$-subdifferential and \eqref{2.11}, we have
\begin{align*}
g(z)-g(z_2) \geq \frac{1}{\tau} \langle D((z_0-D^{-1}v)-z_2-r),z-z_2\rangle- \varepsilon,\forall z \in X,
\end{align*}
i.e.,
\begin{align}
\langle D((z_0-D^{-1}v)-z_2),z_2-z \rangle \geq \tau(g(z_2)-g(z)-\varepsilon)+\langle Dr,z_2-z\rangle,\forall z \in X. \label{2.11a}
\end{align}
Taking $z=z_1$ in \eqref{2.11a}, we have
\begin{align}
\langle D(z_2-z_0+D^{-1}v),z_1-z_2 \rangle\geq \tau(g(z_2)-g(z_1)-\varepsilon)+\langle Dr,z_2-z_1\rangle. \label{2.12}
\end{align}
From \eqref{2.10}, we have
\begin{align}\label{2.12a}
\langle D((z_0-D^{-1}u)-z_1),z_1-z \rangle \geq \tau(g(z_1)-g(z)-\varepsilon).
\end{align}
Setting $z=z_2$ in \eqref{2.12a}, we have
\begin{align}
\langle z_2-z_1,D(z_1-z_0+D^{-1}u) \rangle \geq \tau(g(z_1)-g(z_2)-\varepsilon). \label{2.13}
\end{align}
Adding \eqref{2.12} and \eqref{2.13}, and by a simple manipulation, we get
\begin{align}\label{inq01}
\|z_1-z_2\|_D^2 \leq (\sqrt{2\tau\varepsilon}+\|u-v\|_{D^{-1}})\|z_1-z_2\|_D+2\tau\varepsilon.
\end{align}
Denote
\begin{align}\label{inq02}
t:=\|z_1-z_2\|_D,\,p:=\sqrt{2\tau\varepsilon}+\|u-v\|_{D^{-1}},\,q:=2\tau\varepsilon.
\end{align}
Hence, the inequality $\eqref{inq01}$ is equivalent to the following one:
\begin{align*}
t^2\leq pt+q.
\end{align*}
Thus,
\begin{align}\label{inq03}
(t-\frac{p}{2})^2\leq (\frac{p}{2})^2+q.
\end{align}
Since $t:=\|z_1-z_2\|_D\geq0$, from \eqref{inq03} we have
\begin{align}\label{inq04}
0\leq t\leq\frac{p+\sqrt{p^2+4q}}{2}.
\end{align}
Combining \eqref{inq02} with  \eqref{inq04}, we get
\begin{align*}
0\leq \|z_1-z_2\|_D\leq \frac{1}{2}( \sqrt{2\tau\varepsilon}+\|u-v\|_{D^{-1}}+\sqrt{\|u-v\|_{D^{-1}}^2+10\tau\varepsilon+2\sqrt{2\tau\varepsilon}\|u-v\|_{D^{-1}}}).
\end{align*}
(ii) Note that
\begin{align*}
&\frac{1}{2}\|z-z_0\|_D^2-\frac{1}{2}\|z-z_2\|_D^2\\
&=\frac{1}{2}\|z_2\|_D^2-\frac{1}{2}\|z_0\|_D^2-\langle z_2-z_0,Dz_0\rangle \\
&\quad-\langle z-z_2,v\rangle+\langle z-z_2,D(z_2-z_0+D^{-1}v)\rangle \\
&\geq\frac{1}{2}\|z_2\|_D^2-\frac{1}{2}\|z_0\|_D^2-\langle z_2-z_0,Dz_0\rangle-\langle z_1-z_2,v\rangle \\
&\quad+\langle z_1-z,v\rangle+\tau(g(z_2)-g(z)-\varepsilon)+\langle Dr,z_2-z\rangle,
\end{align*}
where the equality follows from \eqref{equ} and the inequality follows from \eqref{2.11a}. Hence,
\begin{align}  \label{2.15}
\langle z_1-z,v\rangle&\leq \frac{1}{2}\|z-z_0\|_D^2-\frac{1}{2}\|z-z_2\|_D^2-\xi\nonumber\\
 &\quad-\tau(g(z_2)-g(z)-\varepsilon)-\langle Dr,z_2-z\rangle,
\end{align}
where  $\xi=\frac{1}{2}\|z_2\|_D^2-\frac{1}{2}\|z_0\|_D^2-\langle z_2-z_0,Dz_0\rangle-\langle z_1-z_2,v\rangle$.

Also,
\begin{align}\label{2.15a}
\xi&=\frac{1}{2}\|z_2\|_D^2-\frac{1}{2}\|z_0\|_D^2-\langle z_2-z_0,Dz_0\rangle-\langle z_1-z_2,v-u\rangle \nonumber\\
   &\quad+\langle z_2-z_1,D(z_0-z_1)\rangle+\langle z_2-z_1,D(D^{-1}u-z_0+z_1)\rangle                             \nonumber \\
   &\geq\frac{1}{2}\|z_2\|_D^2-\frac{1}{2}\|z_0\|_D^2-\langle z_2-z_0,Dz_0\rangle-\langle z_1-z_2,v-u\rangle \nonumber\\
   &\quad+\langle z_2-z_1,D(z_0-z_1)\rangle+\tau(g(z_1)-g(z_2)-\varepsilon)                                 \nonumber \\
   &=\frac{1}{2}\|z_2-z_1\|_D^2+\frac{1}{2}\|z_0-z_1\|_D^2                             \nonumber \\
   &\quad-\langle z_1-z_2,u-v\rangle+\tau(g(z_1)-g(z_2)-\varepsilon)                             \nonumber   \\
   &\geq \frac{1}{2}\|z_2-z_1\|_D^2+\frac{1}{2}\|z_0-z_1\|_D^2-\frac{1}{2}\|z_2-z_1\|_D^2  \nonumber\\
   &\quad -\frac{1}{2}\|D^{-1}(u-v)\|_D^2+\tau(g(z_1)-g(z_2)-\varepsilon)        \nonumber \\
   &=\frac{1}{2}\|z_0-z_1\|_D^2-\frac{1}{2}\|D^{-1}(u-v)\|_D^2+\tau(g(z_1)-g(z_2)-\varepsilon),
\end{align}
where the first inequality follows from \eqref{2.13} and the second one is due to  the following inequality
\begin{align*}
\langle Dp,q\rangle\leq\frac{1}{2}\|p\|_D+\frac{1}{2}\|q\|_D.
\end{align*}
Combining \eqref{2.15} with \eqref{2.15a}, we have
\begin{align*}
<z_1-z,v>&\leq \frac{1}{2}(\|z-z_0\|_D^2-\|z-z_2\|_D^2)+\frac{1}{2}\|D^{-1}(u-v)\|_D^2 \\
             &\quad-\frac{1}{2}\|z_0-z_1\|_D^2+\tau(g(z)-g(z_1))-\langle Dr,z_2-z\rangle+2\tau\varepsilon \\
             &\leq \frac{1}{2}(\|z-z_0\|_D^2-\|z-z_2\|_D^2)+\frac{1}{2}\|D^{-1}(u-v)\|_D^2 \\
             &\quad-\frac{1}{2}\|z_0-z_1\|_D^2+\tau(g(z)-g(z_1))+\sqrt{2\tau\varepsilon}\|z-z_2\|_D+2\tau\varepsilon,
\end{align*}
where the second inequality follows from the Cauchy-Schwarz inequality.

Therefore, multiplying both sides of the above inequality by $\frac{1}{\tau}$ yields
\begin{align*}
g(z_1)-g(z)+&<z_1-z,\frac{1}{\tau}v>\leq \frac{1}{2\tau}(\|z-z_0\|_D^2-\|z-z_2\|_D^2) \notag \\
&+\frac{1}{2\tau}(\|D^{-1}(u-v)\|_D^2-\|z_0-z_1\|_D^2)+2\varepsilon+\sqrt{\frac{2\varepsilon}{\tau}}\|z-z_2\|_D.
\end{align*}

\end{proof}
\newpage
\section{Main results}\label{sec:alg}
\indent
\setcounter{equation}{0}
In this section, we suppose that $X=\mathbb{R}^m,\,Y=\mathbb{R}^n,\,A\in \mathbb{R}^{n\times m},\,R\in \mathbb{R}^{m\times m}$ and $S\in \mathbb{R}^{n\times n}$. At the same time, we suppose that the matrices $R,\,S$ and $A^TRA$ are symmetric positive definite.

Next we present the inexact primal-dual method for solving \eqref{1.1}.

\begin{algorithm}[h]
\caption{Inexact Primal-Dual Method with Correction Step} 
\textbf{Initialization:} $x^0 \in X,\overline{y}^0 \in Y,\tau_0,\lambda_0>0$.\\
\textbf{Iteration:}\begin{align}
y^{k+1}&\approx_2^{\frac{\varepsilon_{k+1}}{2}} \mathop{\arg\max}_{y\in Y}{L(x^k,y)-\frac{1}{2\tau_k}\|y-\overline{y}^k\|_S^2}           \label{3.1}  \\
x^{k+1}&\approx_2^{\delta_{k+1}}\mathop{\arg\min}_{x\in X}{L(x,y^{k+1})+\frac{1}{2\lambda_k}\|A(x-x^k)\|_R^2}       \label{3.2} \\
\overline{y}^{k+1}&\approx_1^{\frac{\varepsilon_{k+1}}{2}} \mathop{\arg\max}_{y\in Y}{L(x^{k+1},y)-\frac{1}{2\tau_k}\|y-\overline{y}^k\|_S^2} \label{3.3}
\end{align}
\textbf{Until} meet stopping criterion.
\end{algorithm}

Now we consider the two special cases of Algorithm 1.

If we take $\varepsilon_{k+1}=\delta_{k+1}\equiv0$ in Algorithm 1, then Algorithm 1 reduces to the following one:

\begin{algorithm}[h]
\caption{ Primal-Dual Method with Correction Step-A} 
\textbf{Initialization:} $x^0 \in X,\overline{y}^0 \in Y,\tau_0,\lambda_0>0$.\\
\textbf{Iteration:}\begin{align}
y^{k+1}&=\mathop{\arg\max}_{y \in Y}{L(x^k,y)-\frac{1}{2\tau_k}\|y-\overline{y}^k\|_S^2},           \label{3.2a}  \\
x^{k+1}&=\mathop{\arg\min}_{x \in X}{L(x,y^{k+1})+\frac{1}{2\lambda_k}\|A(x-x^k)\|_R^2},      \label{3.2b} \\
\overline{y}^{k+1}&=\mathop{\arg\max}_{y \in Y}{L(x^{k+1},y)-\frac{1}{2\tau_k}\|y-\overline{y}^k\|_S^2}, \label{3.2c}
\end{align}
\textbf{Until} meet stopping criterion.
\end{algorithm}


If, in Algorithm 2, we take $\tau_k=\lambda_k=1$ and
$$R=\left[\begin{matrix}\frac{1}{r_1}I_1 & \textbf{0}  \\
                            \textbf{0}& \frac{1}{r_2}I_2 \end{matrix}\right]
\quad \text{and} \quad \quad
  S=\left[\begin{matrix}\frac{1}{s_1}I_1 & \textbf{0}  \\
                            \textbf{0}& \frac{1}{s_2}I_2\end{matrix}\right]$$
where $r_i,s_i>0,i=1,2$, $I_i(i=1,2)$  are identity matrices, then Algorithm 2 reduces to the following Algorithm 3, which is the PDL method in \cite{Xie2017}.
\begin{algorithm}[h]
\caption{Primal-Dual Method with Correction Step-B} 
\textbf{Initialization:} $x^0 \in X,\overline{y}^0 \in Y,\tau_0,\lambda_0>0$.\\
\textbf{Iteration:}\begin{align}
y^{k+1}&=\mathop{\arg\max}_{y \in Y}{L(x^k,y)-\frac{1}{2}\|y-\overline{y}^k\|_S^2},           \label{3.3a}  \\
x^{k+1}&=\mathop{\arg\min}_{x \in X}{L(x,y^{k+1})+\frac{1}{2}\|A(x-x^k)\|_R^2},      \label{3.3b} \\
\overline{y}^{k+1}&=\mathop{\arg\max}_{y \in Y}{L(x^{k+1},y)-\frac{1}{2}\|y-\overline{y}^k\|_S^2}, \label{3.3c}
\end{align}
\textbf{Until} meet stopping criterion.
\end{algorithm}

\newpage
\begin{remark}\label{re3.2}
It is easy to see that, if $r_is_i<1(i=1,2)$ as required in \cite{Xie2017}, then $R-\overline{\tau}\overline{\lambda}S^{-1}\succ \textbf{0}$ naturally holds. Thus,  our method relaxes the requirement on primal-dual step sizes in \cite{Xie2017}.
\end{remark}
Next we will analyze the convergence of Algorithm 1. Firstly, we prove two important lemmas which will be used in the sequence.
\begin{lemma}\label{le3.1}Let $(y^{k+1},x^{k+1},\overline{y}^{k+1})\in X\times Y\times X$ be obtained from Algorithm 1, then for any pair $(x,y)\in X\times Y$ we have
\begin{align} \label{3.4}
L(x^{k+1},y)-L(x^{k+1},y^{k+1})& \leq \frac{1}{2\tau_k}(\|y-\overline{y}^k\|_S^2-\|y-\overline{y}^k\|_S^2)+\sqrt{\frac{\varepsilon_{k+1}}{\tau_k}}\|y-\overline{y}^{k+1}\|_S \nonumber \\
&\quad+\frac{\tau_k}{2}\|A(x^{k+1}-x^k)\|_{S^{-1}}^2+2\varepsilon_{k+1}- \frac{1}{2\tau_k}\|\overline{y}^k-y^{k+1}\|_S^2.
\end{align}
\end{lemma}
\begin{proof}From Algorithm 1, it is very easy to deduce that the formulas \eqref{3.1} and \eqref{3.2}  are equivalent to the following ones
\begin{align*}
y^{k+1}\approx_2^{\frac{\varepsilon_{k+1}}{2}}\mathop{\arg\min}_{y\in Y}\{g(y)+\frac{1}{2\tau_k}\|y-\overline{y}^k-\tau_kS^{-1}Ax^k\|_S^2\}
\end{align*}
and
\begin{align*}
\overline{y}^{k+1}\approx_1^{\frac{\varepsilon_{k+1}}{2}}\mathop{\arg\min}_{y\in Y}\{g(y)+\frac{1}{2\tau_k}\|y-\overline{y}^k-\tau_kS^{-1}Ax^{k+1}\|_S^2\},
\end{align*}
respectively.

Setting $\tau=\tau_k,\varepsilon=\frac{\varepsilon_{k+1}}{2},D=S,z=y,z_0=\overline{y}^k,z_1=y^{k+1},
z_2=\overline{y}^{k+1},u=-\tau_kAx^k,v=-\tau_kAx^{k+1}$ in Lemma \ref{le2.5}(ii), we get
\begin{align*}
g(y^{k+1})-g(y)+\langle y^{k+1}-y,Ax^{k+1}\rangle &\leq\frac{1}{2\tau_k}(\|y-\overline{y}^k\|_S^2-\|y-\overline{y}^{k+1}\|_S^2)\\
&\quad+\frac{\tau_k}{2}\|A(x^k-x^{k+1})\|_{S^{-1}}^2+2\varepsilon_{k+1}\\
&\quad+\sqrt{\frac{\varepsilon_{k+1}}{\tau_k}}\|y-\overline{y}^{k+1}\|_S-\frac{1}{2\tau_k}\|\overline{y}^k-y^{k+1}\|_S^2.
\end{align*}
Hence,
\begin{align*}
L(x^{k+1},y)-L(x^{k+1},y^{k+1})&=g(y^{k+1})-g(y)-\langle y^{k+1}-y,Ax^{k+1}\rangle\\
&\leq\frac{1}{2\tau_k}(\|y-\overline{y}^k\|_S^2-\|y-\overline{y}^{k+1}\|_S^2) +\frac{\tau_k}{2}\|A(x^k-x^{k+1})\|_{S^{-1}}^2+2\varepsilon_{k+1}\\
&\quad+\sqrt{\frac{\varepsilon_{k+1}}{\tau_k}}\|y-\overline{y}^{k+1}\|_S-\frac{1}{2\tau_k}\|\overline{y}^k-y^{k+1}\|_S^2.
\end{align*}
This completes the proof.
\end{proof}
\begin{lemma}\label{le3.1a}Let $(y^{k+1},x^{k+1},\overline{y}^{k+1})$ be obtained from Algorithm 1, then for any  $x\in X$ we have
\begin{align} \label{3.4a}
L(x^{k+1},y^{k+1})-L(x,y^{k+1})\leq\frac{1}{2\lambda_k}[(\|x-x^k\|_{A^TRA}^2-\|x-x^{k+1}\|_{A^TRA}^2-\|x^{k+1}-x^k\|_{A^TRA}^2]+\delta_{k+1}.
\end{align}
\end{lemma}
\begin{proof}By Definition \ref{d2.3}, the optimal condition of \eqref{3.2} yields
\begin{align*}
\frac{1}{\lambda_k}A^TRA(x^k-x^{k+1}) \in \partial_{\delta_{k+1}}L(x^{k+1},y^{k+1}).
\end{align*}
In view of the definition of $\varepsilon$-subdifferential, we have
\begin{align}\label{3.4b}
L(x^{k+1},y^{k+1})-L(x,y^{k+1})\leq \frac{1}{\lambda_k} \langle A^TRA(x^k-x^{k+1}),x^{k+1}-x)\rangle +\delta_{k+1}.
\end{align}
Setting $a:=x^k, b:=x,c:=x^{k+1},D:=A^TRA$ in \eqref{equ}, we get
\begin{align}\label{3.4c}
\langle A^TRA(x^k-x^{k+1}),x^{k+1}-x)\rangle=-\frac{1}{2}[\|x^k-x^{k+1}\|_{A^TRA}^2+\|x-x^{k+1}\|_{A^TRA}^2-\|x-x^k\|_{A^TRA}^2].
\end{align}
Combining \eqref{3.4b} with \eqref{3.4c}, we know that \eqref{3.4a} holds.
\end{proof}
The following two lemmas play an important role in proving the convergence of Algorithm 1.
\begin{lemma}\label{le3.2}(\cite{SRB3})Assume that the sequence $\{\mu_N\}$ is nonnegative and satisfies the recursion
\begin{align*}
\mu_N^2 \leq T_N+\sum_{n=1}^N \sigma_n \mu_n
\end{align*}
for all $N\geq1$, where $\{T_N\}$ is an increasing sequence, $T_0 \geq \mu_0^2$, and $\sigma_n\geq0$ for all $n\geq0$. Then for all $N\geq1$,
\begin{align*}
\mu_N \leq \frac{1}{2}\sum_{n=1}^N\sigma_n+(T_N+(\frac{1}{2}\sum_{n=1}^N \sigma_n)^2)^{\frac{1}{2}}.
\end{align*}
\end{lemma}

Set
\begin{align} \label{3.7}
\widehat{x}^N:=\frac{1}{N}\sum_{k=0}^{N-1}x^{k+1} \quad \text{and}\quad \widehat{y}^N:=\frac{1}{N}\sum_{k=0}^{N-1}y^{k+1}.
\end{align}
\begin{lemma}\label{th3.1}Let the sequence $\{(x^{k+1},y^{k+1},\overline{y}^{k+1})\}$ be obtained by Algorithm 1 and $(\widehat{x}^N,\widehat{y}^N)$ defined by \eqref{3.7}. Suppose that $\{\tau_k\}$ and $\{\lambda_k\}$ are nondecreasing and $R-\tau_k\lambda_kS^{-1}\succ 0 $. Then for every saddle point  $(x^*,y^*) \in X\times Y$ of \eqref{1.1}, we have
\begin{align}\label{3.9}
L(\widehat{x}^N,y^*)-L(x^*,\widehat{y}^N) \leq \frac{1}{2N\tau_N}[\sqrt{\frac{\tau_N}{\tau_0}}\|y^*-\overline{y}^0\|_S+\sqrt{\frac{\tau_N}{\lambda_0}}\|x^*-x^0\|_{A^TRA}+2A_N+\sqrt{2B_N}]^2,
\end{align}
where  $A_N:=\sum_{k=0}^{N-1} \tau_N \sqrt{\frac{\varepsilon_{k+1}}{\tau_k}}$ and $B_N:=\sum_{k=0}^{N-1} \tau_N(2\varepsilon_{k+1}+\delta_{k+1})$.
\end{lemma}
\begin{proof} Adding \eqref{3.4} and \eqref{3.4a} yields
\begin{align} \label{3.9a}
L(&x^{k+1},y)-L(x,y^{k+1})\leq \frac{1}{2\tau_k}(\|y-\overline{y}^k\|_S^2-\|y-\overline{y}^{k+1}\|_S^2)+\frac{1}{2\lambda_k}
(\|x-x^k\|_{A^TRA}^2-\|x-x^{k+1}\|_{A^TRA}^2)          \notag \\
&+\sqrt{\frac{\varepsilon_{k+1}}{\tau_k}}\|y-\overline{y}^{k+1}\|_S+2\varepsilon_{k+1}+\delta_{k+1}
-\frac{1}{2\lambda_k}\|A(x^k-x^{k+1})\|_{R-\tau_k \lambda_kS^{-1}}^2
-\frac{1}{2\tau_k}\|\overline{y}^k-y^{k+1}\|_S^2.
\end{align}
Since $\{\tau_k\}$ and $\{\lambda_k\}$ are nondecreasing and $R-\tau_k \lambda_kS^{-1} \succ 0$,
\begin{align} \label{3.6}
L(x^{k+1},y)-L(x,y^{k+1})&\leq\frac{1}{2\tau_k}\|y-\overline{y}^k\|_S^2-\frac{1}{2\tau_{k+1}}\|y-\overline{y}^{k+1}\|_S^2 \nonumber \\
&\quad+\frac{1}{2\lambda_k}\|x-x^k\|_{A^TRA}^2-\frac{1}{2\lambda_{k+1}}\|x-x^{k+1}\|_{A^TRA}^2\nonumber\\          &\quad+\sqrt{\frac{\varepsilon_{k+1}}{\tau_k}}\|y-\overline{y}^{k+1}\|_S+2\varepsilon_{k+1}+\delta_{k+1}.
\end{align}
Since $L(x,y)$ and $-L(x,y)$  are convex with respect to $x$ and $y$ respectively, using Jensen inequality and \eqref{3.6}, we have
\begin{align} \label{3.8}
N(L(\widehat{x}^N,y)-L(x,\widehat{y}^N))&\leq\sum_{k=0}^{N-1}L(x^{k+1},y)-L(x,y^{k+1}) \notag \\
&\leq \frac{1}{2\tau_0}\|y-\overline{y}^0\|_S^2-\frac{1}{2\tau_N}\|y-\overline{y}^N\|_S^2+\frac{1}{2\lambda_0}\|x-x^0\|_{A^TRA}^2 \notag \\
&\quad-\frac{1}{2\lambda_N}\|x-x^N\|_{A^TRA}^2          +\sum_{k=0}^{N-1}\sqrt{\frac{\varepsilon_{k+1}}{\tau_k}}\|y-\overline{y}^{k+1}\|_S+\sum_{k=0}^{N-1}(2\varepsilon_{k+1}+\delta_{k+1}).
\end{align}
Setting $x:=x^*$ and $y:=y^*$ in \eqref{3.8} and using \eqref{1.2} we have
\begin{align*}
\|y^*-\overline{y}^N\|_S^2 &\leq \frac{\tau_N}{\tau_0}\|y^*-\overline{y}^0\|_S^2+ \frac{\tau_N}{\lambda_0}\|x^*-x^0\|_{A^TRA}^2 \notag \\
&\quad+\sum_{k=0}^{N-1} \tau_N\sqrt{\frac{\varepsilon_{k+1}}{\tau_k}}\|y^*-\overline{y}^{k+1}\|_S+2\sum_{k=0}^{N-1} \tau_N(2\varepsilon_{k+1}+\delta_{k+1}).
\end{align*}
Set $\mu_N=\|y^*-\overline{y}^N\|_S,T_N=\frac{\tau_N}{\tau_0}\|y^*-\overline{y}^0\|_S^2+
\frac{\tau_N}{\lambda_0}||x^*-x^0||_{A^TRA}^2+2B_N,\sigma_k=2\tau_N\sqrt{\frac{\varepsilon_{k+1}}{\tau_k}}$ in  Lemma \ref{le3.2}.
Obviously,  $T_0 \geq \mu_0^2$ and $\sigma_k \geq 0$. Thus,
\begin{align*}
\|y^*-\overline{y}^N\|_S \leq A_N+(\frac{\tau_N}{\tau_0}\|y^*-\overline{y}^0\|_S^2+
\frac{\tau_N}{\lambda_0}\|x^*-x^0\|_{A^TRA}^2+2B_N+A_N^2)^{\frac{1}{2}}.
\end{align*}
Since $A_N,B_N,\tau_k$ and $\lambda_k$ are nondecreasing,  we have for all $k \leq N$,
\begin{align}\label{3.10}
\|y^*-\overline{y}^k\|_S&\leq A_k+(\frac{\tau_k}{\tau_0}\|y^*-\overline{y}^0\|_S^2+
\frac{\tau_k}{\lambda_0}\|x^*-x^0\|_{A^TRA}^2+2B_k+A_k^2)^{\frac{1}{2}}  \notag \\
&\leq A_N+(\frac{\tau_N}{\tau_0}\|y^*-\overline{y}^0\|_S^2+
\frac{\tau_N}{\lambda_0}\|x^*-x^0\|_{A^TRA}^2+2B_N+A_N^2)^{\frac{1}{2}}  \notag \\
&\leq 2A_N+\sqrt{\frac{\tau_N}{\tau_0}}\|y^*-\overline{y}^0\|_S+
\sqrt{\frac{\tau_N}{\lambda_0}}\|x^*-x^0\|_{A^TRA}+\sqrt{2B_N}.
\end{align}

Hence, setting $x:=x^*$ and $y=:y^*$ in\eqref{3.8},  and using \eqref{3.10} we have
\begin{align*}
N(L(\widehat{x}^N,y^*)-L(x^*,\widehat{y}^N))&\leq \frac{1}{2\tau_0}\|y^*-\overline{y}^0\|_S^2+\frac{1}{2\lambda_0}\|x^*-x^0\|_{A^TRA}^2
+\frac{1}{\tau_N}B_N \\
&\quad+\frac{1}{\tau_N}A_N(2A_N+\sqrt{\frac{\tau_N}{\tau_0}}\|y^*-\overline{y}^0\|_S+
\sqrt{\frac{\tau_N}{\lambda_0}}\|x^*-x^0\|_{A^TRA}+\sqrt{2B_N}) \\
&\leq \frac{1}{2\tau_N}(\sqrt{\frac{\tau_N}{\tau_0}}\|y^*-\overline{y}^0\|_S+\sqrt{\frac{\lambda_N}{\tau_0}}
\|y^*-\overline{y}^0\|_{A^TRA}+2A_N+\sqrt{2B_N})^2,
\end{align*}
which implies  that \eqref{3.9} holds. This  completes the proof.
\end{proof}

\begin{remark}\label{re3.1}
If, in addition,  $A_N$ and $B_N$ are summable and $\{\tau_k\}$ is bounded above, then from \eqref{3.9} we can establish the $O(1/N)$ convergence rate of our method in the ergodic sense.
\end{remark}

\begin{theorem}\label{th3.2} Let $\{x^{k+1},y^{k+1},\overline{y}^{k+1}\}$ be the sequence pair generated by Algorithm 1 and $\{\widehat{x}^N,\widehat{y}^N\}$ be defined by \eqref{3.7} in Lemma \ref{th3.1}. Suppose that the assumptions of Lemma \ref{th3.1} hold and $\tau_k\leq \overline{\tau},\lambda_k\leq \overline{\lambda}$ with $R-\overline{\tau}\overline{\lambda}S^{-1}\succ 0$. If the partial sums $A_N$ and $B_N$ in Lemma \ref{th3.1} are summable and $A$ is of full column rank, then every  cluster point $(\widehat{x},\widehat{y})$ of $\{\widehat{x}^N,\widehat{y}^N\}$ is a saddle point of problem \eqref{1.1}. Moreover, there exists a saddle point $(\widehat{x},\widehat{y}) \in X \times Y$ such that $x^k \rightarrow \widehat{x}$ and $y^k \rightarrow \widehat{y}$ as $k\rightarrow\infty$.
\end{theorem}
\begin{proof}Since $A_N$ and $B_N$ are summable,
\begin{align*}
(\sqrt{\frac{\tau_N}{\tau_0}}&\|y^*-\overline{y}^0\|_S+\sqrt{\frac{\lambda_N}{\tau_0}}
\|y^*-\overline{y}^0\|_{A^TRA}+2A_N+\sqrt{2B_N})^2     \\
\leq (\sqrt{\frac{\overline{\tau}}{\tau_0}}&\|y^*-\overline{y}^0\|_S+\sqrt{\frac{\overline{\lambda}}{\tau_0}}
\|y^*-\overline{y}^0\|_{A^TRA}+2A_N+\sqrt{2B_N} )^2:=C_1<+\infty.
\end{align*}
From \eqref{3.10}, we know that for all $k\leq N$, $\|y^*-\overline{y}^k\|_S\leq C_2<+\infty$. By the same argumentation as for $\overline{y}^k$, from \eqref{3.8} we obtain $\|x^*-x^k\|_{A^TRA}<\infty$ for all $k\leq N$ and hence $\{x^k\}$ is bounded, which implies the boundedness of $\{\widehat{x}^N\}$. Let  $x:=x^*$  and $y=y^*$ in \eqref{3.9a} and then sum the resulting inequality from $k=0$ to $N-1$ to obtain
\begin{align}\label{3.11}
\frac{1}{2\overline{\lambda}} \sum_{k=0}^{N-1}\|x^{k+1}&-x^k\|_{A^T(R-\overline{\tau}\overline{\lambda}S^{-1})A}^2+\sum_{k=0}^{N-1}\frac{1}{2\overline{\tau}}\|\overline{y}^k-y^{k+1}\|_S^2 \leq \frac{1}{2\tau_0}\|y^*-\bar{y}^0\|_S^2-\frac{1}{2\tau_N}\|y^*-\bar{y}^N\|_S^2 \notag \\
&+\frac{1}{2\lambda_0}\|x^* -x^0\|_{A^TRA}^2-\frac{1}{2\lambda_N}\|x^* -x^N\|_{A^TRA}^2
+\frac{C_2}{\tau_N}A_N+\frac{1}{\tau_N}B_N:=C_3<+ \infty.
\end{align}
Letting $N\rightarrow \infty$ in \eqref{3.11} and applying the equivalence of $\|\cdot\|_M$ and $\|\cdot\|_2$, where $M$ denotes the symmetric positive definite matrix, we have
\begin{align} \label{3.12}
x^k-x^{k+1} \rightarrow 0 \quad \text{and} \quad \overline{y}^k-y^{k+1} \rightarrow 0 \quad\textrm{as}\quad k\rightarrow\infty.
\end{align}
Hence, $\{\bar{y}^k-y^{k+1}\}$ is bounded.
Thus,
\begin{align*}
\|y^*-y^{k+1}\|_S\leq \|y^*-\overline{y}^k\|_S+\|\overline{y}^k-y^{k+1}\|_S <+ \infty,
\end{align*}
i.e.,$\{y^k\}$ is bounded, and hence $\{\widehat{y}^N\}$ is also bounded. Hence there exists a subsequence $(\widehat{x}^{N_i},\widehat{y}^{N_i})$ converging to a cluster point$(\widehat{x},\widehat{y})$. Since $f$ and $g$ are l.s.c., from \eqref{3.8} we deduce that, for every fixed $(x,y)\in X\times Y,$
\begin{align*}
L(\widehat{x},y)-L(x,\widehat{y})\leq\liminf_{i \rightarrow \infty} L(\widehat{x}^{N_i},y)-L(x,\widehat{y}^{N_i})=0.
\end{align*}
Taking the supremum over $(x,y)$ in the above inequality  implies that $(\widehat{x},\widehat{y})$ is a saddle point of $L(x,y)$.

Since the sequence pair  $(x^k,y^k)$  is bounded, there exists a subsequence $(x^{k_i},y^{k_i})$  converging to a cluster point $(\widehat{x},\widehat{y})$. Since $(\widehat{x},\widehat{y})$ is a saddle point of $L(x,y)$, by replacing $(x^*,y^*)$ with $(\widehat{x},\widehat{y})$ in \eqref{3.9a}, we know that \eqref{3.12} holds. Hence, $x^{k_i-1}-x^{k_i} \rightarrow 0$ and $\bar{y}^{k_i-1}-y^{k_i} \rightarrow 0$ as $i\rightarrow\infty$. Hence,
\begin{align*}
\|x^{k_i-1}-\widehat{x}\|_{A^TRA} \leq \|x^{k_i-1}-x^{k_i}\|_{A^TRA}+\|x^{k_i}-\widehat{x}\|_{A^TRA} \rightarrow 0\,\,\textrm{as}\,\, i\rightarrow \infty,
\end{align*}
i.e., $x^{k_i-1} \rightarrow \widehat{x}$. Let now $x^{k+1}=H(x^k)$  denote \eqref{3.2b} in Algorithm 2 and
$x^{k+1}=H_{\delta_{k+1}}(x^k)$  denote \eqref{3.2}. In view of the continuity of $H$ , we have
\begin{align*}
\|\widehat{x}-H(\widehat{x})\|_{A^TRA}&=\lim_{i\rightarrow \infty}\|x^{k_i-1}-H(x^{k_i-1})\|_{A^TRA} \\
 &\leq\lim_{i\rightarrow\infty}(\|x^{k_i-1}-H_{\delta_{k_i}}(x^{k_i-1})\|_{A^TRA}
 +\|H_{\delta_{k_i}}(x^{k_i-1})-H(x^{k_i-1})\|_{A^TRA})         \\
 &\leq\lim_{i\rightarrow\infty}(\|x^{k_i-1}-x^{k_i}\|_{A^TRA}+\sqrt{2\bar{\lambda}\delta_{k_i}}) =0,
\end{align*}
where the last inequality follows from  Lemma \ref{le2.3}, Lemma \ref{le2.1} and Definition \ref{2.1}.

By \eqref{3.12} and using the equivalence of the norms $\|\cdot\|_2$ and $\|\cdot\|_S$ we have
\begin{align}\label{eq-a}
\lim_{i\rightarrow\infty}\|\overline{y}^{k_i-1}-y^{k_i}\|_S=0.
\end{align}
In view of  Lemma \ref{le2.5}(i), we have
\begin{align} \label{3.15}
0 \leq \|y^k-\bar{y}^k\|_S& \leq \frac{1}{2}(\sqrt{2\bar{\tau}\varepsilon_k}+
\|x^k-x^{k-1}\|_{A^TS^{-1}A}   \notag\\
&\quad+\sqrt{\|x^k-x^{k-1}\|_{A^TS^{-1}A}^2+10\bar{\tau}\varepsilon_k+
2\sqrt{2\bar{\tau}\varepsilon_k}\|x^k-x^{k-1}\|_{A^TS^{-1}A}}).
\end{align}
Since $x^{k_i-1}-x^{k_i}\rightarrow 0(i\rightarrow\infty)$, taking $k=k_i$ and letting $i\rightarrow\infty$ in the above formula, we get
\begin{align}\label{eq-b}
\lim_{k\rightarrow\infty} \|y^{k_i}-\overline{y}^{k_i}\|_S=0.
\end{align}
Let $y^{k+1}=\Gamma(\overline{y}^k)$ and $\overline{y}^{k+1}=\Psi(\overline{y}^k)$ in Algorithm 1, and
 $y^{k+1}=\Gamma_{\varepsilon_{k+1}}(\overline{y}^k)$ and $\overline{y}^{k+1}=\Psi_{\varepsilon_{k+1}}(\overline{y}^k)$ in Algorithm 2. Thus,
\begin{align}\label{3.13}
\|\widehat{y}-\Gamma \circ \Psi (\widehat{y})\|_S&=\lim_{i\rightarrow\infty}\|y^{k_i}-\Gamma\circ\Psi(\overline{y}^{k_i-1})\|_S\nonumber\\
&\leq \lim_{i\rightarrow\infty}(\|y^{k_i}-\Gamma \circ \Psi_{\varepsilon_{k_i}}(\overline{y}^{k_i-1})\|_S+\|\Gamma \circ \Psi_{\varepsilon_{k_i}}(\overline{y}^{k_i-1})-\Gamma \circ \Psi(y^{k_i-1})\|_S  \nonumber   \\
&\leq \lim_{i\rightarrow\infty}(\|\Gamma_{\varepsilon_{k_i}}(\overline{y}^{k_i-1})-\Gamma (\overline{y}^{k_i})\|_S+\|\Psi_{\varepsilon_{k_i}}(\overline{y}^{k_i-1})- \Psi(\overline{y}^{k_i-1})\|_S      \nonumber\\
&\leq \lim_{i\rightarrow\infty}(\|\Gamma_{\varepsilon_{k_i}}(\overline{y}^{k_i-1})-\Gamma (\overline{y}^{k_i-1})\|_S+\|\Gamma(\overline{y}^{k_i-1})-\Gamma(\overline{y}^{k_i})\|_S
+\sqrt{2\overline{\tau}\varepsilon_{k_i}})    \nonumber  \\
&\leq \lim_{i\rightarrow\infty}[ (2\sqrt{2\overline{\tau}\varepsilon_{k_i}})+ \|\overline{y}^{k_i-1}-\overline{y}^{k_i}\|_S]=\lim_{i\rightarrow\infty} \|\overline{y}^{k_i-1}-\overline{y}^{k_i}\|_S  \nonumber  \\
&\leq \lim_{i\rightarrow\infty} (\|\bar{y}^{k_i-1}-y^{k_i}\|_S+\|\bar{y}^{k_i}-y^{k_i}\|_S)\nonumber \\
&=\lim_{i\rightarrow\infty} \|\overline{y}^{k_i-1}-y^{k_i}\|_S+\lim_{i\rightarrow\infty}\|\overline{y}^{k_i}-y^{k_i}\|_S=0,
\end{align}
where the last equality follows from  \eqref{eq-a} and \eqref{eq-b}.
Hence, from \eqref{3.13} we obtain $\widehat{y}=\Gamma \circ \Psi (\widehat{y})$. Since $x^*=H(x^*)$, it follows that $(\widehat{x},\widehat{y})$ is a fixed point of Algorithm 1 and hence a saddle point of problem \eqref{1.1}. We now use $(x,y)=(\widehat{x},\widehat{y})$ in \eqref{3.6} and sum from $k=k_i,...,N-1$ to obtain
\begin{align}\label{wc}
\frac{1}{2\bar{\lambda}}\|\widehat{x}-x^N\|_{A^TRA}^2+\frac{1}{2\bar{\tau}}\|\widehat{y}-\bar{y}^N\|_S^2
\leq \frac{1}{2\tau_1}\|\widehat{y}-\bar{y}^{k_i}\|_S^2+\frac{1}{2\lambda_1}\|\widehat{x}-x^{k_i}\|_{A^TRA}^2\nonumber\\
+\sum_{k=k_i}^{N-1}\sqrt{\frac{\varepsilon_{k+1}}{\tau_1}}\|\widehat{y}-\bar{y}^{k+1}\|_S+\sum_{k=k_i}^{N-1}(2\varepsilon_{k+1}+\delta_{k+1}).
\end{align}
Since $\varepsilon_k \rightarrow 0$ and $\delta_k\rightarrow 0$ as $k\rightarrow\infty$, the right hand size in \eqref{wc} tends to zero for $i\rightarrow \infty$, which implies that also $x^N \rightarrow \widehat{x}$ and $\bar{y}^N \rightarrow \widehat{y}$ for $N \rightarrow \infty$. Since $x^N \rightarrow \widehat{x}$ as $N\rightarrow\infty$, it is easy to see that $\lim_{N \rightarrow\infty}\|x^N-x^{N-1}\|_{A^TS^{-1}A}=0$. Taking $k=N$ in \eqref{3.15} we have $\|y^N-\bar{y}^N\|_S\rightarrow0$( $N \rightarrow \infty$). Therefore,
\begin{align*}
\lim_{N\rightarrow\infty} \|y^N-\widehat{y}\|_S\leq\lim_{N \rightarrow\infty}(\|y^N-\bar{y}^N\|_S+\|\bar{y}^N-\widehat{y}\|_S)=0.
\end{align*}
Thus, $\|y^N-\widehat{y}\|_S\rightarrow 0$ as $N\rightarrow\infty$ and hence $y^N \rightarrow \widehat{y}$ as $N\rightarrow\infty$. This completes the proof.
\end{proof}
Next we will establish  convergence rates of our method, provided that $\{\delta_k\}$ and $\{\varepsilon_k\}$ decrease like $\mathcal{O}(\frac{1}{k^{2\alpha+1}})$. We first review the following lemma.
\begin{lemma}\label{le3.3}(\cite{RC2020}) For $\omega>-1$, let $s_N:=\sum_{n=1}^N n^\omega$. Then
\begin{align*}
s_N=\mathcal{O}(N^{1+\omega}).
\end{align*}
\end{lemma}
Similar to Corollary 1 in \cite{RC2020} and Theorem 4 in \cite{SSA7}, we have the following theorem.
\begin{theorem}\label{th3.3}Keep the assumptions stated in Theorem \ref{th3.2}. If $\alpha > 0$ and $\delta_k=\mathcal{O}(\frac{1}{k^{2\alpha+1}})$,\,\,$\varepsilon_k=\mathcal{O}(\frac{1}{k^{2\alpha+1}})$, then
\begin{equation*}
L(\widehat{x}^N,y^*)-L(x^*,\widehat{y}^N)=\left \{
      \begin{aligned}
        &\mathcal{O}(\frac{1}{N}),          \quad\quad \,\,\text{$\alpha>\frac{1}{2}$},     \\
        &\mathcal{O}(\frac{\ln^2N}{N}),  \quad \text{$\alpha=\frac{1}{2}$},     \\
        &\mathcal{O}(N^{-{2\alpha}}),    \quad \text{$\alpha\in (0,\frac{1}{2})$}.
        \end{aligned}
\right.
\end{equation*}
\end{theorem}
\begin{proof}If $\alpha>\frac{1}{2}$, then the sequences $\{\varepsilon_k\}$,\,$\{\sqrt{\varepsilon_k}\}$ and $\{\delta_k\}$ are summable and hence $A_N:=\sum_{k=0}^{N-1} \tau_N \sqrt{\frac{\varepsilon_{k+1}}{\tau_k}}$ and $B_N:=\sum_{k=0}^{N-1} \tau_N(2\varepsilon_{k+1}+\delta_{k+1})$ in \eqref{3.9} are bounded. Since the sequence $\{\tau_k\}$ is bounded, from \eqref{3.9} we have
\begin{align*}
L(\widehat{x}^N,y^*)-L(x^*,\widehat{y}^N)=\mathcal{O}(\frac{1}{N}).
\end{align*}
If $\alpha=\frac{1}{2}$, then $\delta_k=\mathcal{O}(\frac{1}{k^2})$ and $\varepsilon_k=\mathcal{O}(\frac{1}{k^2})$. Hence, $B_N$ is bounded. In view of the assumption on $\varepsilon_k$, we have
\begin{align*}
\sqrt{\varepsilon_{k+1}}\leq\frac{c}{k+1}\,\,\, \textrm{for \,\,some}\, \,\,c>0.
\end{align*}
Thus,
\begin{align*}
\sum_{k=0}^{N-1}\sqrt{\varepsilon_{k+1}}\leq\sum_{k=0}^{N-1}\frac{c}{k+1}=c(1+\sum_{k=2}^{N}\frac{1}{k})\leq c(1+\int_1^N\frac{1}{t}dt)=c(1+\ln N).
\end{align*}
Hence, from the boundedness of $\{\tau_k\}$ we know that $A_N=\mathcal{O}(\ln N)$. Therefore, from \eqref{3.9} we get
\begin{align*}
L(\widehat{x}^N,y^*)-L(x^*,\widehat{y}^N)=\mathcal{O}(\frac{\ln^2N}{N}) .
\end{align*}
If $\alpha \in (0,\frac{1}{2})$, then $2\alpha+1>1$ and hence by assumptions $\delta_k=\mathcal{O}(\frac{1}{k^{2\alpha+1}})$ and $\varepsilon_k=\mathcal{O}(\frac{1}{k^{2\alpha+1}})$ we deduce that $B_N$ is bounded. Since $\sqrt{\varepsilon_k}=\mathcal{O}(k^{-\alpha-\frac{1}{2}})$ and $-\alpha-\frac{1}{2}>-1$, from Lemma \ref{le3.3} we obtain $A_N=\mathcal{O}(N^{\frac{1}{2}-\alpha)}$. Thus, from \eqref{3.9} we have
\begin{align*}
L(\widehat{x}^N,y^*)-L(x^*,\widehat{y}^N)=\mathcal{O}(N^{-{2\alpha}}) .
\end{align*}
\end{proof}

\section{Numerical experiments}\label{sec:ex}
\indent
\par
\setcounter{equation}{0}
In this section, we study the numerical solution of  the TV-L$_1$ model for  image deblurring
\begin{align}\label{4.1}
\min_{x \in X}\,\  F(x)=\|Kx-f\|_1 +\mu\|D x\|_1,
\end{align}
where $f\in Y$ is a given (noisy) image, $K: X\rightarrow Y$ is a known linear (blurring) operator, $D:X\rightarrow Y$ denotes the gradient operator and $\mu$ is a regularization parameter. Now we introduce the variables $\gamma_1,\gamma_2 $, which satisfy $\gamma_1,\gamma_2>0$ and $\gamma_1+\gamma_2=\mu$. Then, \eqref{4.1} can be written as
\begin{align*}
\min_{x \in X} \|Kx-f\|_1 +\gamma_1\|D x\|_1+\gamma_2\|D x\|_1.
\end{align*}

Further, the above formula can be rewritten as \cite{CP2016}
\begin{align}\label{4.2}
\min_{x \in X} \max_{y \in Y} L(x,y):=\gamma_1\|D x\|_1+\langle Ax,y\rangle-\delta_{C_1}(u)-\delta_{C_{\gamma_2}}(v)-\langle f,u\rangle,
\end{align}
where $C_{\lambda}=\{y\in Y|\, \|y\|_{\infty}\leq\lambda \}$, $y=\left [\begin{matrix}u  \\v \end{matrix} \right]$, and
$A=\left [\begin{matrix}K  \\ \gamma_2 D \end{matrix} \right]$.
Next we suppose that $\mathcal{N}(K)\bigcap\mathcal{N}(D) =\{0\}$, where $\mathcal{N}(Q)$ represents the null space of the matrix $Q$, and this assumption has been used in many references including similar types of problems, for example, in \cite{WYYZ2008}. Under this assumption,  $A$ is of full column rank.

For simplicity we set  $\tau_k=\lambda_k=1$, $R=\left[\begin{matrix}\frac{1}{r_1}I_1 & \textbf{0}  \\
                           \textbf{0}& \frac{1}{r_2}I_2 \end{matrix}\right]
\quad \text{and} \quad \quad
 S=\left[\begin{matrix}\frac{1}{s_1}I_1 & \textbf{0}  \\
                           \textbf{0}& \frac{1}{s_2}I_2\end{matrix}\right]$ in Algorithm 1. Setting $\varepsilon_{k+1}=0$, and we can compute $y^{k+1}$ by the following formula:
\begin{align*}
u^{k+1}&=P_{C_1}(\bar{u}^k+s_1(Kx^k-f)), \\
v^{k+1}&=P_{C_{\gamma_2}}(\bar{v}^k+s_2\gamma_2 Dx^k),
\end{align*}
where $y^{k+1}=(u^{k+1},v^{k+1})^T$. Similarly,  we can  get $\bar{y}^{k+1}=(\overline{u}^{k+1},\overline{v}^{k+1})$ where $\overline{u}^{k+1}$ and $\overline{v}^{k+1}$ can be obtained  by replacing $x^k$ with $x^{k+1}$ in $u^{k+1}$ and $v^{k+1}$, respectively.

Now we consider the computation of the following subproblem:
\begin{align}\label{4.3}
x^{k+1}&\approx_2^{\delta_{k+1}}\arg\min_{x\in X}{L(x,y^{k+1})+\frac{1}{2}\|A(x-x^k)\|_R^2}\nonumber\\
&\approx_2^{\delta_{k+1}}\arg\min_{x\in X}{\gamma_1\|D x\|_1+\langle Ax,y^{k+1}\rangle+\frac{1}{2}\|A(x-x^k)\|_R^2}.
\end{align}

The above formula can be equivalently rewritten as
\begin{align}\label{4.3a}
x^{k+1} \approx_2^{\delta_{k+1}}\arg\min_{x \in X}{\gamma_1\|D x\|_1+\frac{1}{2}\|Bx-\xi\|^2},
\end{align}
where $B=\left [\begin{matrix}\frac{1}{\sqrt{r_1}}K  \\\frac{\gamma_2}{\sqrt{r_2}}D \end{matrix} \right]$ and $\xi=\left [\begin{matrix} \frac{1}{\sqrt{r_1}}Kx^k- \sqrt{r_1}u^{k+1}\\\frac{\gamma_2}{\sqrt{r_2}}Dx^k-\sqrt{r_2}v^{k+1} \end{matrix} \right]$. We note that $B^TB=\frac{1}{r_1}K^TK+\frac{\gamma_2^2}{r_2}D^TD$ is symmetrically positive definite because $B$ is of full column rank. Therefore, there exists $z \in X$ such that $Bz=\xi$. Further, \eqref{4.3} can be rewritten equivalently as
\begin{align}\label{4.3a}
x^{k+1} \approx_2^{\frac{\delta_{k+1}}{\gamma_1}}\arg\min_{x \in X}{\|D x\|_1+\frac{1}{2\gamma_1}\|x-z\|^2_{B^TB}}.
\end{align}

Next we will show that the subproblem \eqref{4.3a} can be computed by approximately minimizing duality gap.

Setting $H(x):=h(Dx)=\|D x\|_1,\varphi(x)=\frac{1}{2\gamma_1}\|x-z\|^2_{B^TB}$,  we consider the following primal problem:
\begin{align}\label{4.4a}
\min_{x \in X}h(Dx)+\varphi(x):=\|D x\|_1+\frac{1}{2\gamma_1}\|x-z\|^2_{B^TB}
\end{align}
together with its dual problem:
\begin{align} \label{4.4}
\min_{v \in Y} \varphi^*(-D^T v)+h^*(v)=\frac{1}{2\gamma_1}
\|\gamma_1D^Tv-B^T\xi\|^2_{(B^TB)^{-1}}
-\frac{1}{2\gamma_1}\|\xi\|_2^2+\delta_\Omega(v)
\end{align}
which $\Omega=\{v | \|v\|_{\infty}\leq 1\}$. If $\bar{v}$ is a solution of Problem \eqref{4.4}, then
\begin{align}\label{4.4b}
\bar{x}=(B^TB)^{-1}(B^T\xi-\gamma_1D^T\bar{v})
\end{align}
is a solution of Problem \eqref{4.4a}. Since $h$ is positively homogeneous, from Remark 1 of \cite{VSBV2013} we get $h^*(v)=H^*(D^Tv)$. Now we  consider the dual gap
\begin{align}\label{4.5}
\Psi(\bar{x},\bar{v})&=h(D\bar{x})+\varphi(\bar{x})
+\varphi^*(-D^T \bar{v})+h^*(\bar{v})\nonumber\\
&=H(\bar{x})+H^*(D^T\bar{v})+\langle-\bar{x},D^T\bar{v}\rangle\nonumber\\
&=\sup_{v \in Y} \{\langle v,\bar{x}\rangle-H^*(v)\}-\langle\bar{x},D^T\bar{v}\rangle+H^*(D^T\bar{v})\nonumber\\
&\geq H^*(D^T\bar{v})-H^*(v)+ \langle v-D^T\bar{v},\bar{x}\rangle,\,\,\forall\,\,v\in Y.
\end{align}
 Thus, if $\Psi(\bar{x},\bar{v}) \leq \delta$, where $\delta>0$ is  some given tolerance, then $\bar{x} \in \partial_{\delta} H^*(D^T\bar{v}) $, which by Theorem 2.4.4 of \cite{Zal2002} is equivalent to $D^T\bar{v} \in \partial_{\delta} H(\bar{x})$, and hence from Definition \ref{d2.3} this implies $\bar{x}\approx_2^{\delta} \textrm{Prox}_{\gamma_1 H}^{B^TB}(z)$. Therefore,
\begin{align}\label{4.5a}
\Psi(\bar{x},\bar{v}) \leq \delta\Rightarrow \bar{x}\approx_2^{\delta} \textrm{Prox}_{\gamma_1 H}^{B^TB}(z),
\end{align}
where $\bar{x}$ and $\bar{v}$ satisfy \eqref{4.4b}.
  Thus, we use FISTA method (\cite{BT2009}) to solve the dual problem \eqref{4.4} so as to better evaluate
the gap. In view of \eqref{4.5a}, we adopt the following inequality as the stopping criterion of inner loop:
 \begin{align}\label{4.6}
 \Psi(x^{k+1},v^{k+1}) \leq \delta_{k+1},
 \end{align}
 where  $\delta_{k+1}=O(1/(k+1)^{2\alpha+1}(\alpha>0)$.  In the following we will report the numerical experiment results.

The MATLAB codes are run on a PC (with CPU Intel i5-5200U) under MATLAB Version 8.5.0.197613 (R2015a) Service Pack 1. We report numerical results of the proposed methods. We test the images cameraman.png$(256 \times 256)$ and man.png($1024\times 1024$), as presented in Figures 1 and 2. At the same time, we adopt the following stopping rule:
\begin{align*}
\frac{F(x^k)-F(x^*)}{F(x^*)}<10^{-5},
\end{align*}
where $x^*$ is a solution of  the TV-L$_1$ model \eqref{4.1}.

\begin{figure}[h]
\begin{minipage}[t]{0.5\linewidth}
\centering
\includegraphics[width=2.2in]{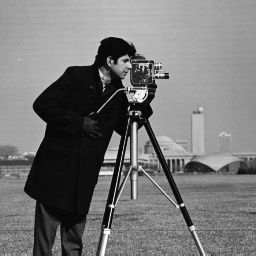}
\caption{Original cammeraman.png(256$\times$256)}
\label{fig:side:a}
\end{minipage}%
\begin{minipage}[t]{0.5\linewidth}
\centering
\includegraphics[width=2.2in]{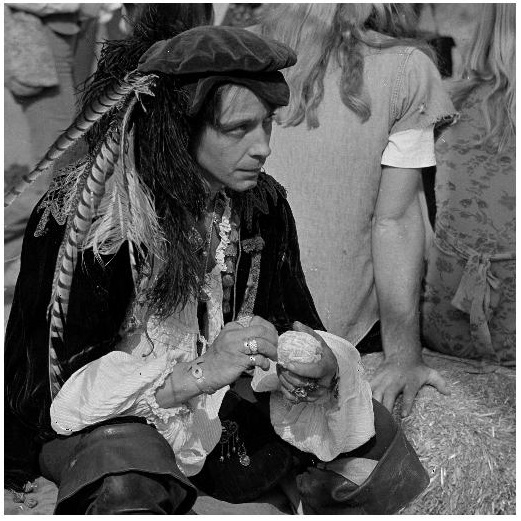}
\caption{Original man.png(1024$\times$1024)}
\label{fig:side:b}
\end{minipage}
\end{figure}

\subsection{Sensitivity of parameters}

In this section, we will analyze the sensitivity of parameters.  In this test, average blur with hsize=9 was applied to  the original image cameraman.png(see Figure 3) (256$\times$256) by fspecial(’average’,9), and 20\% salt-pepper noise was added in. According to Theorem \ref{th3.3}, the convergence \textcolor{red} {rate} of Algorithm 1 depends on the value of parameter $\alpha$. At the same time,  from \eqref{4.6} we know that the iteration number of inner loop closely relates to the parameter $\alpha$.  Hence, we first study the sensitivity of $\alpha$. In the following  experiment, we take $\mu=0.05,\,s_1=1,\,s_2=2,\,r_1=\frac{0.99}{s_1},\,r_2=\frac{0.99}{s_2}$, $\gamma_1=\gamma_2=\frac{1}{2}\mu$. We choose the $256\times256$ cammeraman.png as the test picture and take $\alpha \in (0,2.3)$. The iteration number of outer loop is fixed as $100$. In Figure 5, the ordinate denotes the iteration number of the inner loop while  the abscissa denotes the value of $\alpha$. From Figure 5, we can see that, when $\alpha$ is not very large, the iteration number of inner loop is very little; However, as $\alpha$ increases, the iteration number of inner loop also increases rapidly. Similar results can also be found in \cite{RC2020}.
\begin{figure}[h]
\begin{minipage}[t]{0.5\linewidth}
\centering
\includegraphics[width=2.2in]{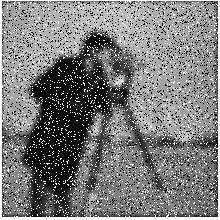}
\caption{Cammeraman.png with noise}
\label{fig:side:a}
\end{minipage}%
\begin{minipage}[t]{0.5\linewidth}
\centering
\includegraphics[width=2.2in]{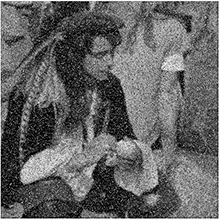}
\caption{Man.png with noise}
\label{fig:side:b}
\end{minipage}
\end{figure}

\begin{figure}[H]
\centering \includegraphics[width=.5\textwidth]{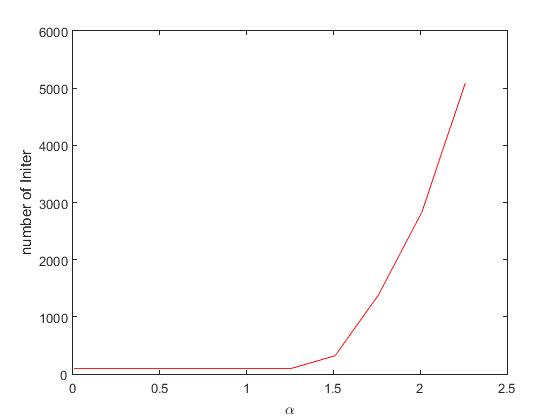}
\caption{Sensitivity of $\alpha$}
\label{fig1}
\end{figure}

Next we investigate the sensitivity of parameters $s_1$ and $s_2$. We still take $\mu=0.05,\,r_1=\frac{0.99}{s_1},\,r_2=\frac{0.99}{s_2},\,\gamma_1=\gamma_2=\frac{1}{2}\mu$ and fix $\alpha=1$. If $s_2$ is fixed as 2, we take $s_1\in (0.8,2.5)$; If $s_1$ is fixed as 1, we take $s_2\in (0.8,2.3)$. The iteration number of outer loop is fixed as $100$. In Figures 6 and 7, the ordinate denotes the running time of Algorithm 1 while  the abscissa denotes the value of $s_1$ or $s_2$. From Figures 6 and 7, we can see that,  the running time  decreases as $s_1$ or $s_2$ increases.

\begin{figure}[H]
\begin{minipage}[t]{0.5\linewidth}
\centering
\includegraphics[width=2.2in]{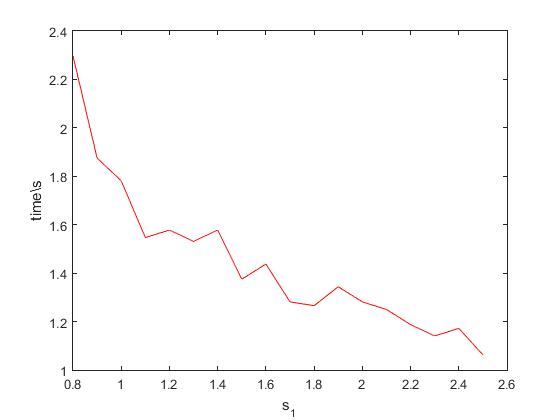}
\caption{Sensitivity of $s_1$}
\label{fig:side:a}
\end{minipage}%
\begin{minipage}[t]{0.5\linewidth}
\centering
\includegraphics[width=2.2in]{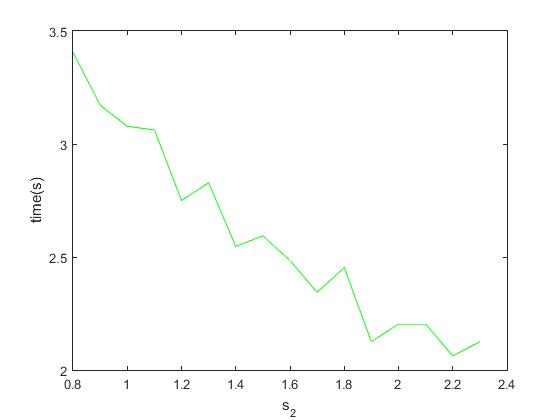}
\caption{Sensitivity of $s_2$}
\label{fig:side:b}
\end{minipage}
\end{figure}

Further, we consider the sensitivity of parameters $\gamma_1$ and $\gamma_2$ which satisfy $\gamma_1+\gamma_2=\mu$. We take $\alpha=1$, $\mu=0.05$ and fix the iteration number of outer loop as $200$.  In Figure 8, the ordinate denotes the value of $F(x^k)-F(x^*)$ over $F(x^*)$ while the abscissa denotes the value after the number of inner iterations is taken $\textrm{log}_{10}$. In addition, "initr"  denotes the number of total iterations of inner loop(i.e., the second subproblem).   At the same time, we choose $\gamma_1$ as four different values $\frac{1}{2}\mu, \frac{1}{2.5}\mu, \frac{1}{3}\mu$ and $\frac{1}{3.5}\mu$  which correspond to four different curves in Figure 8, respectively.  These curves indicate that, when the value of $\gamma_1$ or $\gamma_2$ varies, the number of inner iterations  changes remarkably. Hence, the CPU time increases  markedly as $\gamma_1$ decreases. By testing Figure 3, the CPU time corresponding to the above four choices of $\gamma_1$ is $95.8s,149.7s,210.9s$ and $279.2s$, respectively.
\begin{figure}[htbp]
\centering \includegraphics[width=.5\textwidth]{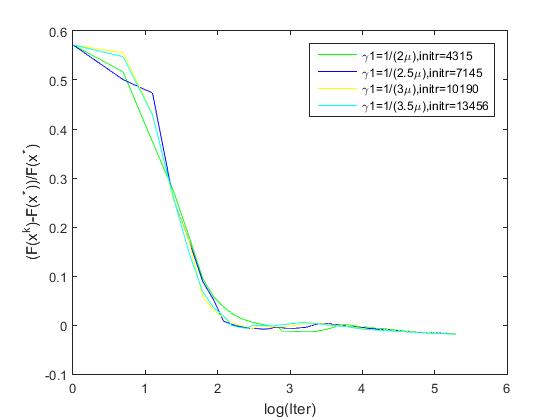}
\caption{Sensitivity of $\gamma$}
\label{fig1}
\end{figure}

Finally, we analyze the variation of Algorithm 1's numerical performance
with respect to various choices of  $\alpha$ when $\gamma_1$ is fixed as  $\frac{1}{3}\mu$. Besides, average blur with hsize=9 was applied to  the original image man.png(see Figure 4) (1024$\times$1024) by fspecial(’average’,9), and 20\% salt-pepper noise was added in. In Table 1,  “CPU”, "iter-out" and "iter-in" denote the CPU time in seconds, the iteraion number of outer and inner loops, respectively. Testing Figures 3 and 4 yields the following results:

\qquad\qquad\qquad Table 1
\begin{table}[!htbp]
\centering
\begin{tabular}{|c|ccc|ccc|}
\hline
 &&Figure 3&&&Figure 4&\\
\hline
$\alpha$ &CPU&iter-out &iter-in &time(s)&iter-out&iter-in\\
\hline
0.1   &0.8438&10&18&26.0900&10&18  \\
\hline
0.3   &1.0625&10&34&35.2031&10&32  \\
\hline
0.5   &1.2500&9&35&40.4688&10&41  \\
\hline
0.8   &1.5625&10&51&48.3594&10&50  \\
\hline
1   &1.8281&11&63&58.3906&11&63  \\
\hline
\end{tabular}
\end{table}
\\
From Table 1 we can see that, as $\alpha$ increases,  the CPU time and the number of inner iterations increase while the number of outer iterations keeps invariant bascially.

\subsection{Image denoising}
In this section, we apply our Algorithm 1 to  image debluring of TV-L$_1$ model \eqref{4.1}. In the following tables and figures, "CP", “iCP”, "PDL" and "iPDL"  denote Algorithm 1 in \cite{CP2011}, Algorithm (4.2) in \cite{RC2020}, Algorithm PDL in \cite{Xie2017} and our Algorithm 1, respectively.   In this experiment, we test Figures 3 and 4. We fixed the number of iterations as 200 and the penalty coefficient $\mu=0.1$. When the above four Algorithms  are implemented, their respective parameters are given in Table 2.

\qquad\qquad\qquad Table 2
\begin{table}[H]
\centering
\begin{tabular}{|c|c|}
\hline
$CP$    &$\tau=\sigma=\frac{0.99}{\sqrt{8}}$\\
\hline
$iCP$   &$\tau=\sigma=0.99,\alpha=1 $  \\
\hline
$PDL$   &$s_1=2,r_1=\frac{0.99}{s_1},s_2=1,r_2=\frac{0.99}{s_2}$\\
\hline
$iPDL$  &$s_1=2,r_1=\frac{0.99}{s_1},s_2=1,r_2=\frac{0.99}{s_2},\alpha=1,\gamma_1=\gamma_2=\frac{1}{2\mu}$\\
\hline

\end{tabular}
\end{table}

The  restored images by the above four Algorithms are displayed in Figure 9. Obviously, our algorithm and CP algorithm get better  restoration quality compared with the iCP and PDL methods. In our experiment,  we find that, if we increase the number of iterations to 1000 or more, all four algorithms can restore the image  with almost the same quality, but our algorithm needs fewer iterations than other three ones.

\begin{figure}[h]
\centering \includegraphics[width=.9\textwidth]{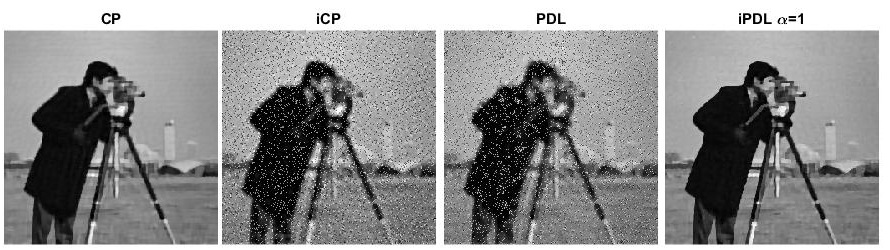}
\caption{Restored images}
\label{fig1}
\end{figure}

\section{Conclusions}\label{sec:con}
In this paper, we propose an inexact primal-dual method for the saddle point problem by applying  inexact extended proximal operators. We show the convergence of our Algorithm 1, provided that the partial sums $A_N$ and $B_N$  are summable. The $O(1/N)$ convergence rate in the ergodic sense is also established. We also apply our method to solve TV-L$_1$ image deblurring  problems and verify their efficiency numerically.

It is worth mentioning that our method has some existing algorithms as special cases by the appropriate choices of parameters. Besides,  our method also relaxes the requirement on primal-dual step sizes, for example, in \cite{Xie2017}. At present, however, we are not able to  provide the  accelerated versions of our method , for example, under the assumption that $f$ or $g$ is strongly convex. Hence, this will be the subject of future research.


\begin{thebibliography}{5}

\bibitem{BT2009}A. Beck,and M. Teboulle, A fast iterative shrinkage-thresholding algorithm for linear inverse problems. SIAM J Imaging Sciences,  2(1)(2009), 183-202.
\bibitem{BPCPE2010}S. Boyd , N. Parikh , E. Chu , B. Peleato, and J. Eckstein, Distributed Optimization and
Statistical Learning via the Alternating Direction Method of Multipliers[J], Foundations and
Trends in Machine Learning, 3(2010), 1-122.
\bibitem{CCS2010} J.-F. Cai, E.J. Candès, and Z. Shen, A singular value thresholding algorithm for matrix completion.
SIAM Journal on Optimization, 20(4) (2010), 1956-1982.
\bibitem{CP2011}A. Chambolle, and T. Pock. A first-order primal-dual algorithm for convex problems with applications
to imaging. Journal of Mathematical Imaging and Vision, 40(1)(2011), 120-145.
\bibitem{CP2016}A. Chambolle, and T. Pock. An introduction to continuous optimization for imaging. Acta Numerica,
25(2016), 161-319.
\bibitem{CR2013}Patrick L. Combettes, and Noli N. Reyes, Moreau's decomposition in Banach spaces, Mathematical Programming(Ser. B),  139(2013), 103-114.
\bibitem{EB1992}J. Eckstein, and D. P. Bertsekas, On the Douglas-Rachford splitting method and the proximal
point algorithm for maximal monotone operators, Mathematical Programming, 55 (1992), 293-318.
\bibitem{EB2016}M. Ehrhardt, M.Betcke,  Multicontrast MRI reconstruction with structure-guided total variation.
SIAM Journal on Imaging Sciences, 9(3)(2016), 1084-1106.
\bibitem{EZC2010} E. Esser, X. Zhang, and T. F. Chan. A general framework for a class of first order primal-dual
algorithms for convex optimization in imaging science. SIAM Journal on Imaging Sciences, 3(4)(2010), 1015-1046.
\bibitem{FP2011}J. M. Fadili, and G. Peyre,  Total variation projection with frst order schemes. IEEE Trans. Image Process. 20(3)(2011), 657-669.
\bibitem{HTY2012} B.S. He, M. Tao, and X. Yuan, Alternating direction method with Gaussian back substitution for separable convex programming, SIAM Journal on Optimization, 22(2)(2012), 313-340.
\bibitem{HY2012}B.S. He, X. Yuan, Convergence analysis of primal-dual algorithms for a saddle-point problem: from contraction perspective, SIAM Journal on Imaging Sciences, 5(1)(2012), 119-149.
\bibitem{HYY2014}B.S. He, Y.F. You, and X.M. Yuan, On the convergence of primal-dual hybrid gradient
algorithm, SIAM Journal on Imaging Sciences, 7 (2014), 2526-2537.
\bibitem{LXY2020}Y. Liu, Y. Xu and W. Yin, Acceleration of primal-dual methods by preconditioning and simple subproblem procedures, Journal of Scientific Computing, 86, 21 (2021). https://doi.org/10.1007/s10915-020-01371-1
\bibitem{MGC2011} S. Ma,  D. Goldfarb, L. Chen,  Fixed point and Bregman iterative methods for matrix rank minimization. Mathematical Programming, 128(1)(2011), 321-353.
\bibitem{NWY2011}M. K. Ng, F. Wang, and X. Yuan, Inexact alternating direction methods for image recovery,
SIAM Journal on Scientific Computing, 33 (2011), 1643-1668.
\bibitem{PCBC2009}T. Pock, D. Cremers, H. Bischof, and A. Chambolle, An algorithm for minimizing the
Mumford-Shah functional, in Computer Vision, 2009 IEEE 12th International Conference
on, IEEE, 2009, 1133-1140.
\bibitem{RC2020}J. Rasch J, and A. Chambolle,  Inexact first-order primal-dual algorithms. Computational Optimization and Applications,  76(2020), 381-430.
\bibitem{ROF1992}L. Rudin, S. Osher, and E. Fatemi, Nonlinear total variation based noise removal algorithms, Phys.D, 60 (1992), 227-238.
\bibitem{SSA7}S. Salzo, and S. Villa. Inexact and accelerated proximal point algorithms[J]. Journal of Convex Analysis,19(4)(2012), 1167-1192.
\bibitem{SRB3} M.Schmidt, N.L. Roux, and F.R. Bach,  Convergence rates of inexact proximal-gradient methods for
convex optimization. In: Shawe-Taylor, J., Zemel, R.S., Bartlett, P.L., Pereira, F., Weinberger, K.Q.
(eds.) Advances in Neural Information Processing Systems, vol. 24, pp. 1458-1466. Curran Associates Inc, Red Hook (2011)
\bibitem{Val2014}T. Valkonen, A primal–dual hybrid gradient method for nonlinear operators with applications
to mri, Inverse Problems, 30 (2014), 055012.
\bibitem{VSBV2013}S. Villa, S. Salzo, L. Baldassarre, and A. Verri. Accelerated and inexact forward-backward algorithms.SIAM Journal on Optimization,   23(3)(2013), 1607-1633.
\bibitem{WYYZ2008}Y. Wang, J. Yang, W. Yin, and Y. Zhang, A new alternating minimization algorithm for total variation image reconstruction. SIAM Journal on Imaging Sciences, 1(3)(2008), 248-272.
\bibitem{Xie2017}Z. Xie, A primal-dual method with linear mapping for a saddle point problem in image deblurring. J. Vis. Commun. Image R., 42(2017), 112-120.
\bibitem{Zal2002}C. Zalinescu. Convex Analysis in General Vector Spaces. World Scientific, 2002
\bibitem{ZC2008}M. Zhu, and T. Chan, An efficient primal-dual hybrid gradient algorithm for total variation
image restoration, UCLA CAM Report, 34 (2008).



\end{thebibliography}
\end{document}